\documentclass[10pt]{amsart}
\usepackage{latexsym, amsmath,amssymb}

\usepackage{xcolor}

\renewcommand{\epsilon}{\varepsilon}
\newcommand{\newsection}[1]
{\subsection{#1}\setcounter{theorem}{0} \setcounter{equation}{0}
\par\noindent}

\newtheorem{theorem}{Theorem}

\newtheorem{lemma}[theorem]{Lemma}
\newtheorem{corr}[theorem]{Corollary}

\newtheorem{proposition}[theorem]{Proposition}
\newtheorem{deff}[theorem]{Definition}

\newcommand{\bth}{\begin{theorem}}
\newcommand{\ble}{\begin{lemma}}
\newcommand{\bcor}{\begin{corr}}

\newcommand{\bdeff}{\begin{deff}}

\newcommand{\bprop}{\begin{proposition}}
\newcommand{\ele}{\end{lemma}}
\newcommand{\ecor}{\end{corr}}
\newcommand{\edeff}{\end{deff}}

\newcommand{\eprop}{\end{proposition}}

\newcommand{\Rn}{{\mathbb R}^n}

\newcommand{\la}{\lambda}

\renewcommand{\Pi}{\varPi}

\renewcommand{\epsilon}{\varepsilon}

\newcommand{\one}{{\bf 1}}

\begin{document}

\title
{Curvature and harmonic analysis on compact manifolds}

\thanks{The  author was supported in part by the NSF (DMS-1665373). }


\author{Christopher D. Sogge}
\address{Department of Mathematics,  Johns Hopkins University,
Baltimore, MD 21218}
\email{sogge@jhu.edu}

\begin{abstract}
We discuss problems that relate curvature and concentration properties of eigenfunctions and quasimodes
on compact boundaryless Riemannian manifolds.  
These include new  sharp $L^q$-estimates, $q\in (2,q_c]$, 
$q_c=2(n+1)/(n-1)$,  of log-quasimodes that characterize compact
connected space forms in terms of the growth rate of  $L^q$-norms of  such quasimode for these relatively small Lebesgue exponents $q$.
 No such characterization is possible for  any  exponent $q> q_c$.
\end{abstract}

\maketitle

\newsection{Introduction and main results}

Since at least the 1970s there has been a great deal
of interest in exploring the role of curvature in problems naturally arising in harmonic analysis.  Prominent among these was the Stein-Tomas Fourier
restriction theorem \cite{TomasWill}.  This result says that if $n\ge2$ there is a uniform constant $C=C_n$ so that if $d\omega$ denotes surface
measure on the unit sphere
\begin{equation}\label{1}
\|f\|_{L^2(S^{n-1})}=\bigl(\, \int_{S^{n-1}}|\Hat f(\omega)|^2 \, d\omega \, \bigr)^{1/2}
\le C\|f\|_{L^p(\Rn)}, \, \, 1\le p\le \tfrac{2(n+1)}{n+3}, \, \, f\in {\mathcal S}(\Rn).
\end{equation}
This uniform inequality and a simple density argument allow one to define the restriction of the Fourier transform of $f\in L^p(\Rn)$ to
$S^{n-1}$ as an element of $L^2(S^{n-1})$ if $1\le p\le \tfrac{2(n+1)}{n+3}$.  This  is remarkable since the Fourier transform of
elements of $f\in L^p(\Rn)$ for $p>1$ are only defined almost everywhere, and  the range of exponents in \eqref{1} is sharp due to the
Knapp example.  It is the curvature of $S^{n-1}$ that allows \eqref{1}.  Indeed, if one replaced $S^{n-1}$ by any hypersurface containing
a nontrivial open subset of a hyperplane, then it is an easy exercise to see that \eqref{1} can never hold for {\em any} exponent $p>1$.

In \cite{sogge88} the author initiated the program of attempting to obtain natural generalizations of \eqref{1} in the setting of compact
boundaryless Riemannian manifolds $(M,g)$ of dimension $n\ge2$.  Specifically, if $\Delta_g$ is the associated Laplace-Beltrami operator
and $P=\sqrt{-\Delta_g}$ it was shown in \cite{sogge88} that one has the universal spectral projection bounds
\begin{multline}\label{2}
\bigl\| \one_{[\la,\la+1]}(P)f\|_{L^q(M)}\lesssim \la^{\mu(q)} \|f\|_{L^2(M)}, \, \, \la\ge 2, 
\\
\text{with } \quad \mu(q)=
\begin{cases}
n(\tfrac12-\tfrac1q)-\tfrac12, \, \, q\ge q_c=\tfrac{2(n+1)}{n-1},
\\
\tfrac{n-1}2(\tfrac12-\tfrac1q), \, \, 2\le q\le q_c.
\end{cases}
\end{multline}
Here $\one_{[\la,\la+1]}(s)$, $s\in {\mathbb R}$, is the indicator function of the unit-length interval $[\la,\la+1]$ and 
$\one_{[\la,\la+1]}(P)$ is the operator defined by the spectral theorem.  One calls $q_c=\tfrac{2(n+1)}{n-1}$ the ``critical exponent''
for \eqref{2} since the estimates for $q\in (2,\infty] \backslash \{q_c\}$ easily follow from the special case of \eqref{2} where
$q=q_c$.  Also, one can see that the bounds in \eqref{2} are sharp for {\em any} compact manifold $M$ (see \cite{SFIO2}).

In the Euclidean case the spectral projection operators associated with the intervals $[\la, \la+\delta]$ are just
$$\bigl(\one_{[\la,\la+\delta]}(\sqrt{-\Delta_{\Rn}}) f\bigr)(x)=(2\pi)^{-n}\int_{|\xi|\in [\la,\la+\delta]} e^{ix\cdot \xi} \Hat f(\xi) \, d\xi,$$
and, by an easy argument which makes use of duality and scaling, one sees that \eqref{1} is equivalent to the  uniform bounds
for $\la\ge2$ and $\delta\in (0,1]$
\begin{equation}\label{3}
\bigl\|\one_{[\la,\la+\delta]}(\sqrt{-\Delta_{\Rn}})f\|_{L^q(\Rn)}
\lesssim \la^{\mu(q)} \cdot 
\begin{cases}
\delta^{1/2} \|f\|_{L^2(\Rn)}, \, \, q\in [q_c,\infty],
\\
\delta^{\frac{n+1}2(\frac12-\frac1q)}\|f\|_{L^2(\Rn)}, \, \, q\in (2,q_c].
\end{cases}
\end{equation}
It is also easy to see that the bounds in \eqref{3} are optimal, and if $\delta=1$ these bounds agree with those in \eqref{2}.

As we mentioned, the unit-band spectral projection bounds \eqref{2} for compact manifolds are always sharp.  Notwithstanding, one can
ask whether, if, as in \eqref{3},  the unit intervals $[\la,\la+1]$ are replaced by smaller ones, say, $[\la,\la+\delta]$ with $\delta=\delta(\la)
\searrow 0$,  is it possible to improve upon the bounds in \eqref{2}?  On the sphere, there can be no improvement as was shown by the author
in \cite{sogge86}.  On the other hand, Zelditch and the author \cite{SoggeZelditchMaximal} showed that for generic manifolds there are
improved estimates for $q>q_c$, and, arguably introduced the program of attempting to find geometries for which improvements of
\eqref{2} are possible, not just for this range of exponents, but any $q>2$.

Implicit in B\'erard~\cite{Berard} is that if $\delta(\la)=(\log\la)^{-1}$ then 
$\|\one_{[\la,\la+\delta(\la)]}\|_{2\to \infty}=
(\delta(\la))^{1/2}\la^{\frac{n-1}2}$, $\la\ge2$, if the sectional curvatures of $M$ are all nonpositive, which agrees with the
$\delta^{1/2}$ improvements  in the first part of \eqref{3}.  In a later important paper, Hassell and Tacy~\cite{HassellTacy}
extended this result to all ``supercrtical'' exponents showing that, under this curvature assumption, one  has
\begin{equation}\label{4}
\|\one_{[\la,\la+\delta(\la)]}\|_{2\to q}\le C_q (\delta(\la))^{1/2}\la^{\mu(q)}, \quad \text{if  } \, \,
q>q_c, \, \, \text{and } \, \, \delta(\la)=(\log\la)^{-1}.
\end{equation}
The choice of  $\delta=\delta(\la)=(\log\la)^{-1}$ naturally arises due to the role of the Ehrenfest time in manifolds of negative sectional 
curvature (cf. Zelditch~\cite{Zelquantum}).

It is an easy exercise involving the Cauchy-Schwarz inequality to show that, in \eqref{4} (as in \eqref{3}), the dependence on $\delta(\la)$ is sharp.  On the
other hand, until recently, it was not known whether improvements of \eqref{2} can hold for $q\in (2,q_c]$.  The ``critical'' exponent 
$q=q_c$ seemed especially difficult to handle since one has to rule out quasimodes behaving like the ``zonal functions'' and 
``Gaussian beams'' on $S^n$ under, say, curvature assumptions for $(M,g)$.  Recall (cf. \cite{sogge86}) that zonal functions
saturate the bounds for $S^n$ in \eqref{2} for $q\ge q_c$ and the Gaussian beams saturate the bounds for $q\in (2,q_c]$.  Both saturate
the bounds for $q=q_c$ and have very different ``profiles''.  The zonal functions are highly concentrated near the poles on $S^n$, while
the Gaussian beams are highly concentrated near the equator, which, of course is a periodic geodesic.  So, any attempts to improve
upon the bounds in \eqref{2}, and try to extend those in \eqref{4} to other exponents, necessarily would seem to have to rule out both types of concentration
for $q=q_c$ and concentration near geodesics for the range $q\in (2,q_c)$.

In a joint paper with Zelditch \cite{SoggeZelditchL4} we were able to obtain for $n=2$ partial, but non-optimal improvements of \eqref{2} for intervals
$[\la,\la+\delta(\la)]$ with $\delta(\la)\searrow 0$
when $q\in (2,q_c)$.  After this, Blair and the author, in a series papers \cite{BlairSoggeRefined}, \cite{blair2015refined}
and \cite{BSTop}, extended such results to higher dimensions.  These results failed to address the
more difficult problem of obtaining improved spectral projection bounds for shrinking spectral windows when $q=q_c$.

In \cite{sogge2015improved}, the author showed that, if all the sectional curvatures of $M$ are nonpositive, then there are
very weak ``log-log'' improvements of bounds for the spectral projection operators in \eqref{4} if $q=q_c$.  These bounds
were then improved considerably by  Blair and the author in \cite{SBLog},
 showing that, under this curvature assumption,  one has the bounds
\begin{equation}\label{5}
\|\one_{[\la,\la+\delta(\la)]}\|_{2\to q_c}\le C (\delta(\la))^{\sigma_n}\la^{\mu(q_c)}, \quad \text{if  } \, \,
 \, \delta(\la)=(\log\la)^{-1}.
\end{equation}
The powers $\sigma_n>0$ that were obtained were not optimal and went to zero as $n\to \infty$; however, they were the first logarithmic
improvements to be obtained, which allowed a partial extension of the results \eqref{4} of B\'erard~\cite{Berard} and Hassell and Tacy~\cite{HassellTacy}.  The proof of \eqref{5} was
simplified considerably, and stronger results were recently obtained by Blair, Huang and the author \cite{BHSsp}.

In ongoing work we are able to further simplify the arguments and finally obtain the following {\em sharp} estimates.

\begin{theorem}\label{thm1.1}  Let $(M,g)$ be an $n$-dimensional connected compact Riemannian manifold.
Then, if all the sectional curvatures are {\em nonpositive}, for $\la\gg1$ we have the uniform bounds
\begin{equation}\label{1.7}
\bigl\| \chi_{[\la, \la+(\log\la)^{-1}]} f\bigr\|_{L^q(M)}\le C\bigl(\la (\log\la)^{-1}\bigr)^{\mu(q)} \|f\|_{L^2(M)}, \, \, \,
2<q\le q_c,
\end{equation}
with $q_c$ and $\mu(q)$ as in \eqref{2}.  Moreover, if all the sectional
curvatures of $M$ are {\em negative}, for $\la\gg 1$ we have the uniform bounds
\begin{equation}\label{1.8}
\bigl\| \chi_{[\la, \la+(\log\la)^{-1}]} f\bigr\|_{L^q(M)}\le C_q \, \la^{\mu(q)} (\log\la)^{-1/2} \|f\|_{L^2(M)}, \, \, \,
2<q\le q_c,
\end{equation}
with the constant $C_q$ in \eqref{1.8} depending on $q$.
\end{theorem}

Note that the bounds in \eqref{1.8} are stronger than the corresponding sharp Euclidean ones in \eqref{3} if
$2<q<q_c$.   The fact that such estimates are valid on compact manifolds but not $\Rn$ is surprising.  On the other hand
Chen and Hassell \cite{ChenHassell} obtained analogous bounds in their extension of the Stein-Tomas 
restriction theorem to ${\mathbb H}^n$ (see also \cite{SHuangSogge} for the first type of such extension dealing with $q\ge q_c$).

The sharpness of the bounds is due to the following result which also tells us that it is possible to characterize compact
space forms in terms of the growth rate of  $L^q$-norms of log-quasimodes for relatively small exponents $q$.  To state it we let $V_{[\la,\la+(\log\la)^{-1}]}$ denote
those functions whose $P$--spectrum  lies in $[\la,\la+(\log\la)^{-1}]$.  This is the
space of all log-quasimodes associated to a given frequency $\la \ge2$.
We also recall that if $f,g\ge0$ then 
$f(\la)=\Theta(g(\la))$ if $\limsup_{\la\to \infty}\tfrac{f(\la)}{g(\la)}\in (0,\infty)$, i.e., $f(\la)=O(g(\la))$ and also
$f(\la)=\Omega(g(\la))$ (the negation of $f(\la)=o(g(\la)$).  

Our other main result then is the following result which
says that compact space forms are characterized by the growth rate of $L^q$-norms of
$L^2$-normalized log-quasimodes,
if  $q\in (2,q_c]$.

\begin{theorem}\label{thm1.2}  Assume that $(M,g)$ is a connected compact manifold of constant
sectional curvature $K$ and fix {\em any} exponent
$q\in (2,q_c]$.
Then, if $\mu(q)$ is as in \eqref{2},
\small
\begin{equation}\label{shape}
\sup_{\Phi_\la\in V_{[\la,\la+(\log)^{-1}]}}\frac{\|\Phi_\la\|_{L^q(M)}}{\|\Phi_\la\|_{L^2(M)}} =
\begin{cases}
\Theta(\la^{\mu(q)}(\log\la)^{-1/2}) \, \, \iff \, \, K<0
\\
\Theta(\la^{\mu(q)}(\log\la)^{-\mu(q)}) \, \, \iff \, \, K=0
\\
\Theta(\la^{\mu(q)}) \, \, \iff \, \, K>0.
\end{cases}
\end{equation}
\normalsize
Also, if $(\log\la)^{-1}\le \delta(\la)\searrow 0$ as $\la\to \infty$ and $\la\to \la \, \delta(\la)$ is non-decreasing for $\la\ge2$,
\small
\begin{equation}\label{shape2}
\sup_{\Phi_\la\in V_{[\la,\la+\delta(\la)]}}\frac{\|\Phi_\la\|_{L^q(M)}}{\|\Phi_\la\|_{L^2(M)}} =
\begin{cases}
\Theta(\la^{\mu(q)}(\delta(\la))^{1/2}) \, \, \iff \, \, K<0
\\
\Theta(\la^{\mu(q)}(\delta(\la))^{\mu(q)}) \, \, \iff \, \, K=0
\\
\Theta(\la^{\mu(q)}) \, \, \iff \, \, K>0.
\end{cases}
\end{equation}
\normalsize
\end{theorem}

We should point out that the  estimates \eqref{4} of Hassell and Tacy~\cite{HassellTacy}, which are optimal in terms of their $\delta(\la)$--dependence,  do not distinguish between
compact manifolds of zero sectional curvatures from ones all of whose sectional curvatures are negative.  So, interestingly,
one must use the range of exponents $q\in (2,q_c]$ as in Theorem~\ref{thm1.2} to characterize compact space forms in terms of the growth
rate of log-quasimodes.  Also, the harmonic analysis that is used to treat the two cases of $q>q_c$ or $q\in (2,q_c]$ is much
different.  The latter requires adaptations of  bilinear techniques from Lee~\cite{LeeBilinear} and Tao, Vargas and Vega~\cite{TaoVargasVega}.

Note that the bounds in Theorem~\ref{thm1.1} along with the universal bounds in \eqref{2} say that one automatically
has the analog of \eqref{shape} with $\Theta(\, \cdot \, )$ replaced by $O(\, \cdot \, )$.  To prove that one also has
the $\Omega(\, \cdot \, )$ lower bounds and thus obtain the preceding theorem one needs to use the constant 
curvature assumptions.  Proving these lower bounds when $K<0$ or $K>0$ is straightforward; however, handling
flat space forms is more difficult.  For the $K=0$ case one needs to see that for any flat compact manifolds one can
obtain a Knapp type example using ideas from Brooks~\cite{BrooksQM} and Zelditch and the author~\cite{SoggeZelditchL4}.

\newsection{Applications and further problems}

We can use the above bounds to obtain new results concerning concentration properties of log-quasimodes in various geometries.
One is improvements of the universal lower bounds for $L^1$-norms of the author and Zelditch~\cite{SoggeZelNodal} that 
were used to obtain (non-optimal) lower bounds for the size of nodal sets of eigenfunctions on compact manifolds.  The universal lower
bounds say that
$\la^{-\frac{n-1}4}\lesssim \|\Phi_\la\|_{L^1(M)}$, if
$\|\Phi_\la\|_2=1$  and the spectrum of $\Phi_\la$ lies in
$ [\la,\la+1]$.
The Gaussian beams on $S^n$ saturate these lower bounds.
Using Theorem~\ref{thm1.1} one can improve these lower bounds as follows
\begin{equation}\label{2.1}
\small
\|\Phi_\la\|_{L^1(M)}  \gtrsim \la^{-\frac{n-1}4} \cdot 
\begin{cases}
(\log\la)^{\frac{n-1}4}, \,  \text{if all the sectional curvatures are nonpositive}
\\
(\log\la)^{N}, \, \forall \, N, \, \text{if all the sectional curvatures are negative}.
\end{cases}
\normalsize
\end{equation}
provided that  $\text{Spec }\Phi_\la \subset [\la,\la+(\log\la)^{-1}]$ and $\|\Phi_\la\|_2=1$.

It is straightforward to see that the lower bounds in \eqref{2.1} are always sharp if $M$ is flat.  An interesting but potentially difficult
problem would be to see to what extent one could improve the lower bounds for manifolds of negative curvature, either for the above
quasimodes or even eigenfunctions with eigenvalue $\lambda$.

Another interesting problem would be to to see whether the bounds in \eqref{1.8} are valid under weaker assumptions, such, as for
instance, the assumption that the geodesic flow is Anosov.  It would also be  interesting to see whether there are improvements
of the universal bounds in \eqref{2} for shrinking spectral widows for {\em generic} manifolds for $q\in (2,q_c)$ or the more difficult
case where $q=q_c$.  This would complement the results in \cite{SoggeZelditchMaximal}.

\bibliography{refs}
\bibliographystyle{abbrv}

\end{document}